\documentclass{article}
\usepackage{amsmath,amsthm,amsfonts}
\usepackage{enumerate}
\usepackage{amssymb}
\usepackage{epsf,psfig}

\newcommand{\seconds}{{\mathrm{s}}}
\newcommand{\miliseconds}{{\mathrm{ms}}}
\newcommand{\hertz}{{\mathrm{Hz}}}

\title{A simple mechanism for balancing at the border of instability with
applications to persistent neural activity}

\author{Luc Moreau\\[1ex]
Dynamics and Control group\\ Mechanical Engineering department\\ Eindhoven University of Technology\\ The Netherlands\\[2ex]
Eduardo Sontag\\[1ex]
Department of Mathematics\\ Rutgers, The State University
of New Jersey\\ Piscataway~N.J.\ 08854, USA}

\date{\today}

\begin{document}

\maketitle

\begin{abstract}

Some biological systems operate at the critical point between
stability and instability and this requires a fine-tuning of
parameters.  We bring together two examples from the literature that
illustrate this: neural integration in the nervous system and hair
cell oscillations in the auditory system.  In both examples the
question arises as to how the required fine-tuning may be achieved and
maintained in a robust and reliable way.  We study this question using
tools from nonlinear and adaptive control theory.  We illustrate our
approach on a simple model which captures some of the essential
features of neural integration. As a result, we propose a large class
of feedback adaptation rules that may be responsible for the
experimentally observed robustness of neural integration.  We mention
extensions of our approach to the case of hair cell oscillations in
the ear.
\end{abstract}

\vspace*{1cm}
Persistent neural activity is prevalent throughout the nervous system.
Numerous experiments have demonstrated that persistent neural activity
is correlated with short-term memory.  A
prominent example concerns
the oculomotor system---see \cite{Ro:89,AkGaSeBaTa:01} for a review and experimental facts.
The brain moves the eyes with quick saccadic movements. Between
saccades, it keeps the eyes still by generating a continuous and
constant contraction of the eye muscles; thus requiring a constant
level of neural activity in the motor neurons controlling the eye
muscles.  This constant neural activity level serves as a short-term
memory for the desired eye position.  During a saccade, a brief burst
of neural activity in premotor command neurons induces a persistent
change in the neural activity of the motor neurons, via a mechanism
equivalent to integration in the sense of calculus.  Neural activity
of an individual neuron, however, has a natural tendency to decay with
a relaxation time of the order of milliseconds.  Therefore the
question arises as to how a transient stimulus can cause persistent
changes in neural activity.  According to a long-standing hypothesis,
persistent neural activity is maintained by synaptic feedback loops.
Positive feedback can oppose the tendency of a pattern of neural
activity to decay.  If the feedback is weak, then the natural tendency
to decay dominates and neural activity decays. As the feedback
strength is increased, the neural dynamics undergo a bifurcation and
become unstable.  When the feedback is tuned to exactly balance the
decay, then neural activity neither increases nor decreases but
persists without change.  This, however, requires a fine-tuning of the
synaptic feedback strength and the question arises as to how a
biological system can achieve and maintain this
fine-tuning~\cite{ArRo:91,Se:96,SeLeReTa:neuron00,SeLeReTa:jcn00}.
Some gradient descent and function approximation algorithms performing
this fine-tuning have been proposed~\cite{ArRo:91,SeLeReTa:neuron00}
and a feedback learning mechanism based on differential anti-Hebbian
synaptic plasticity has been studied in~\cite{XiSe:00}.  Nevertheless,
it is still unclear how the required fine-tuning is physiologically
feasible. For this reason, a different model for neural integration
based upon bistability has recently been proposed
in~\cite{KoRaKeLi:02}. In the present paper, we do not follow the line
of research based upon bistability. Instead, we pursue the hypothesis
of precisely tuned synaptic feedback. The present paper proposes an
adaptation mechanism that may be responsible for the fine-tuning of
neural integrators and that may explain the experimentally observed
robustness of neural integrators with respect to perturbations. Before
we present this adaptation mechanism in detail, we first discuss a
similar phenomenon in the auditory system.

In order to detect the sounds of the outside world, hair cells in the
cochlea operate as nanosensors which transform acoustic stimuli into
electric signals.  In \cite{CaDuJuPr:00,EgOsChHuMa:00,Ju:01,JuAnDu:01}
these hair cells are described as active systems capable of generating
spontaneous oscillations.  Ions such as $\mathrm{Ca}^{++}$ are
believed to contribute to the hair cell's tendency to self-oscillate.
For low concentrations of the ions, damping forces dominate and the
hair cell oscillations are damped. As the concentration increases the
system undergoes a Hopf bifurcation, the dynamics become unstable, and
the hair cells exhibit spontaneous oscillations. In
\cite{CaDuJuPr:00,EgOsChHuMa:00,Ju:01,JuAnDu:01} the hair cells are
postulated to operate near the critical point, where the activity of
the ions exactly compensates for the damping effects.  As before, this
requires a fine-tuning of parameters (the ion concentrations) and
again the question arises as to how this fine-tuning can be achieved
and maintained.  In
\cite{CaDuJuPr:00,Ju:01} a feedback mechanism has been proposed which
could be responsible for maintaining this fine-tuning.

It thus seems that operating in the vicinity of a bifurcation is a
recurrent theme in biology. And the question as to how proximity to
the bifurcation point may be achieved and maintained in a noisy
environment may be of considerable, general interest. We view the
two presented examples as special instances of the following general
problem.  Consider a forced dynamical system, described by a
differential equation $\dot{x}=f_\mu(x,u(t))$.  The right-hand side of
this equation depends on a parameter~$\mu$ and the unforced dynamics
$\dot{x}=f_\mu(x,0)$ are assumed to exhibit a bifurcation when $\mu$
equals a critical value~$\mu_0$.  The problem consists of finding a
feedback adaptation rule for the parameter~$\mu$ which guarantees
proximity to the bifurcation point; that is, which steers $\mu$
toward its critical value~$\mu_0$.  This adaptation law may depend on
$\mu$ and $x$ but should be independent of $\mu_0$, since this
critical value is not known precisely.  This abstract
formulation captures common
features of both biological examples and suggests some
unexpected links with the literature.  Questions very similar to the
present one have been studied extensively in the literature on
adaptive control~\cite{KrKaKo:95} and stabilization~\cite{So:98}; and
the general problem is closely related
to extremum seeking~\cite{KrWa:00}, and to
instability detection~\cite{AnShFlKeRiKr:99}, where an
operating parameter is adapted on-line in order to experimentally locate
bifurcations.

Although the above general formulation is convenient,
there is little
hope that a complete and satisfactory theory can be developed that
applies to all possible instances of the problem.
Simplifying assumptions make it
more tractable. In this letter, we study in detail what is
probably the most simple but nontrivial instance of the general
problem.  We consider the one-dimensional system
\begin{equation}
\dot{x}=- \mu_0 x + \mu x + u(t),\qquad(\dot{x}=\mathrm{d}x/\mathrm{d}t)
\label{e:ni}
\end{equation}
which captures some of the essential features of neural integration
and is in fact closely related to the autapse model from
\cite{SeLeReTa:jcn00}. With this interpretation, $x$ is a strictly
positive variable representing neural activity in the integrator
network and $u(t)$ represents the signal generated by the premotor
command neurons.  The term $-\mu_0 x$ corresponds to the natural decay
of neural activity and $\mu x$ represents a positive, synaptic
feedback loop.  Of course, when studying neural integration, questions
can be investigated at varying levels of detail. It is clear that a
simple model as (\ref{e:ni}) has several limitations.  Because of its
one-dimensional nature, the present model is, for example, unable of
reproducing the distributed nature of persistent activity patterns
observed in the brain.  Nevertheless, Eq.~(\ref{e:ni})
captures a key feature of neural integration: when the feedback is
tuned to exactly balance the decay, Eq.~(\ref{e:ni}) behaves as an
integrator and produces persistent neural activity.
Eq.~(\ref{e:ni}) is therefore a valuable model when studying
fine-tuning of neural integrator networks~\cite{Ro:89,Sh:89}.

We are interested in the fine-tuning of Eq.~(\ref{e:ni}) and
study this question using tools from nonlinear and adaptive control
theory.  First, we ignore the presence of the input~$u(t)$
and consider the simpler equation
\begin{equation}
\dot{x}=- \mu_0 x + \mu x.
\label{e:pa}
\end{equation}
We present a large class of feedback adaptation laws for~(\ref{e:pa})
which steer $\mu$ to its critical value~$\mu_0$; thus enabling the
automatic self-tuning of parameters and the spontaneous generation
of persistent neural activity.  We consider adaptation
laws
%
\footnote{
The dynamics for $\mu$ could come from synaptic plasticity.  In
particular, the term $f(x)$ might be related to types of synaptic
plasticity that depend on the temporal ordering of presynaptic and
postsynaptic spiking, as in \cite{XiSe:00}.
}
%
of the form
\begin{equation}
\dot{\mu}=f(x)-g(\mu).
\label{e:adaptation}
\end{equation}
We show that, under three very mild conditions, this adaptation rule
guarantees convergence to the bifurcation point for~(\ref{e:pa}).  The
first condition requires that~$g$ is a strictly increasing function.
This means that the term~$-g(\mu)$ in~(\ref{e:adaptation}) acts as a
negative feedback. As a consequence, if the neural activity~$x$ were
constant in~(\ref{e:adaptation}), then the synaptic feedback
gain~$\mu$ would naturally relax to a rest value depending on~$x$
via the equation~$f(x)=g(\mu)$. The second condition states that there
exists~$x^*$ such that $f(x^*)=g(\mu_0)$. This condition implies that,
if the neural activity would be constant and equal to~$x^*$
in~(\ref{e:adaptation}), then the synaptic feedback gain~$\mu$ would
naturally relax to its critical, desired value~$\mu_0$.  Of course
there is no guarantee that the neural activity would be equal to, or
even converge to, this special value~$x^*$. Instead, the level of
neural activity is governed by Eq.~(\ref{e:pa}). Therefore, in order
for the adaptation law~(\ref{e:adaptation}) to work, we need to impose
a last condition, that~$f$ is a decreasing
function. This means that the level of neural activity negatively
regulates the synaptic feedback strength.

We now show that, under these three conditions, the feedback
adaptation law~(\ref{e:adaptation}) indeed tunes the synaptic feedback
gain~$\mu$ to exactly balance the natural decay rate~$\mu_0$. We begin
with noticing that the combined system of
equations~(\ref{e:pa})--(\ref{e:adaptation}) has a unique rest
point. This equilibrium is determined by setting the right-hand sides
of~(\ref{e:pa})--(\ref{e:adaptation}) equal to zero, yielding $x=x^*$
and $\mu=\mu_0$.  Although the precise value of~$\mu_0$ is unknown,
if we are able to prove that all trajectories
of~(\ref{e:pa})--(\ref{e:adaptation}) converge to this (unknown) fixed
point, then it follows that~$\mu$ indeed converges to its desired,
critical value~$\mu_0$. In order to prove this, we introduce a
coordinate transformation~$q=\ln(x)-\ln(x^*)$ and~$p=\mu-\mu_0$. This
transforms~(\ref{e:pa})--(\ref{e:adaptation}) into $\dot{q}=p$ and
$\dot{p}=f(\exp(q)x^*)-g(p+\mu_0)$. In these new coordinates, the
dynamics take the form of a nonlinear mass-spring-damper system (with
unit mass, nonlinear spring characteristic~$f(\exp(\cdot)x^*)$ and
nonlinear damping function~$g(\cdot+\mu_0)$).  It follows from
physical energy considerations that this system exhibits damped
oscillations~\cite{MoSoMu:arxiv}. This shows that all trajectories
of~(\ref{e:pa})--(\ref{e:adaptation}) indeed converge to the unique
fixed point, where $\mu=\mu_0$.

The above coordinate
transformation reveals a subtle relationship between self-tuning of
bifurcations and the internal model principle (``IMP'') from robust control
theory (see~\cite{YiHuSiDo:00,So:imp} for a discussion of the IMP
from a systems biology perspective).  This relation is
made explicit by the equation~$\dot{q}=p$, which represents an
integrator and corresponds to integral action studied in robust
control theory.
One regards the constant~$\mu_0$ as an
unknown perturbation acting on the system.
The IMP implies that, in order to track this constant perturbation, the
system dynamics should contain integral action.
The integral action is generated by the biological system
itself, and not by the feedback adaptation law.

We have so far ignored the presence of the
signal~$u(t)$. We showed that the adaptation
law~(\ref{e:adaptation}) tunes the synaptic feedback gain to exactly
compensate for the natural decay rate, resulting in the spontaneous
generation of persistent neural activity. At these equilibrium
conditions, the action potential firing rate equals~$x^*$, which is
related to~$\mu_0$ by $f(x^*)=g(\mu_0)$.  In the next paragraphs, we
take into account the effect of the input~$u(t)$. In this case, the
value~$x^*$ will play the role of a parameter that influences the
accuracy with which the feedback adaptation law guarantees proximity
to the bifurcation point.

The signal~$u(t)$ will in general result in a time-varying action
potential firing rate~$x(t)$. The mechanism with which this happens,
is determined by the neural integrator equation~(\ref{e:ni}) and the
adaptation law~(\ref{e:adaptation}).  For the purpose of analysis, we
make two simplifying assumptions, both of which seem to be natural and
physically relevant for neural integration.  First, we assume that,
over any sufficiently large time interval~$[t_0,\,t_0+T]$, the time
spent by~$x(t)$ in any interval~$[x_1,\,x_2]$ is approximately
independent of~$t_0$.  In more mathematical terms, we assume the
existence of a function~$P(x)$ such that for every test
function~$\alpha(x)$, the time average
\(
\frac{1}{T}\int_{t_0}^{t_0+T}\alpha(x(t))\,\mathrm{d}t
\)
converges to
\(
\int_0^\infty P(x)\alpha(x)\,\mathrm{d}x
\)
as $T\rightarrow\infty$, uniformly with respect to~$t_0$.  Secondly,
we assume that the adaptation law acts on a much slower time scale
than the time variations in~$x(t)$.  Under these assumptions, the
effect of the action potential firing rate~$x(t)$ on the adaptation
law~(\ref{e:adaptation}) may be approximated by the average effect
$\dot{\mu}=\int_{0}^{\infty}P(x)f(x)\,\mathrm{d}x-g(\mu)$.  It is now
clear when the adaptation law guarantees proximity to the bifurcation
point:  if the {compatibility
condition}~$\int_{0}^{\infty}P(x)f(x)\,\mathrm{d}x=f(x^*)$ is
satisfied, then time scale separation arguments suggest that $\mu$ will
converge approximately to $\mu_0$ and the neural integrator will
approximately behave as a perfect integrator.
The compatibility condition may by interpreted as
follows~\cite{MoSoMu:arxiv}.  When the premotor command signal $u(t)$
has zero time-average and the adaptation law acts on a slow time
scale, then Eq.~(\ref{e:ni}) behaves as a good integrator and the
firing rate $x(t)$ equals the time-integral of $u(t)$ plus an
integration constant. The compatibility condition ensures that this
integration constant is compatible with the desired range for the
firing rate $x(t)$.

\begin{figure}[h]
\hspace*{\fill}
\psfig{file=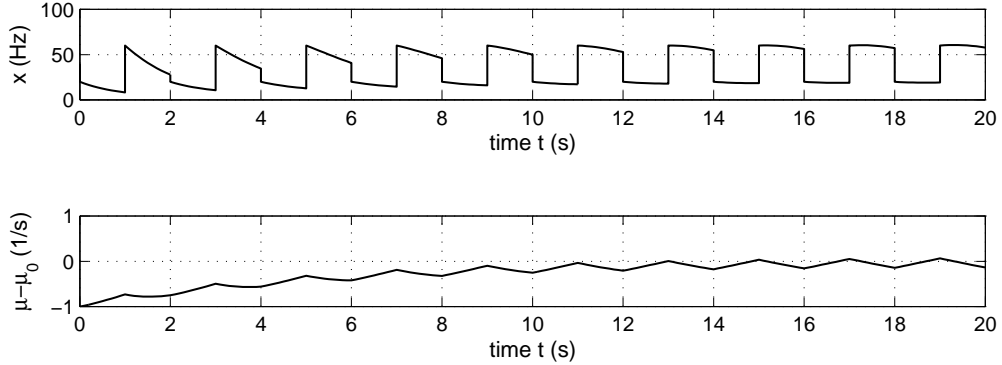,height=5cm}
\hspace*{\fill}
\caption{\label{f2} Tuning of a neural integrator. Simulation of
Eq.~(\ref{e:ni}) with $\dot{\mu}=\varepsilon(-ax-b\mu+c)$ and
$u(t)=\sum_{t_i}\delta(t-t_i)\{x_\mathrm{desired}(t_i+)-x(t_i-)\}$,
where $\delta(\cdot)$ represents the Dirac impulse and where the sum
goes over all saccade times~$t_i$. The constants are:
$\mu_0=200~\seconds^{-1}$, $a=1~\seconds^{-1}$,
$b=0.01~\seconds^{-1}$, $c=42~\seconds^{-2}$ and
$\varepsilon=0.01$. The adaptation law satisfies the compatibility
condition~$a(20~\hertz+60~\hertz)/2+b\mu_0=c$.}
\end{figure}

We illustrate this result on a particular example representative for
saccadic eye movements.  We consider the case of periodic saccadic eye
movements asking for an action potential firing rate in the motor
neurons alternating between $20~\hertz$ and $60~\hertz$ every second.
At each saccade, a brief burst of neural activity in premotor command
neurons changes the actual firing rate. We assume that this change is
such that immediately after each saccade, the actual firing rate
equals the desired firing rate.  Between saccades, we assume that no
input is applied
%
\footnote{
We assume that no feedback is applied to keep $x$ at its desired level
between saccades, which is consistent with experimental observations.
}.
%
If the neural integrator is perfectly tuned, then
the actual firing rate will remain constant between saccades and equal
to the desired firing rate (eyes are fixed). If the neural integrator
is not perfectly tuned, then the actual firing rate will deviate from
the desired firing rate (eyes drift) until a new saccade occurs which
brings the actual firing rate to its new desired value.
%
%
Fig.~\ref{f2} shows the
results of a simulation where the adaptation law satisfies the
compatibility condition of the previous paragraph.
In the beginning of the simulation, we have mis-tuned the neural
integrator.  Clearly, after a short transient, the adaptation law
achieves excellent tuning and the drift between two successive
saccades becomes negligible.

We have thus shown that an adaptation law can tune a neural integrator
with great accuracy to its bifurcation point. In order to achieve
perfect tuning, however, the adaptation law itself needs to satisfy a
compatibility condition.~%
%
%
It seems that we
have merely moved the problem of fine-tuning from the neural
integrator to the adaptation law.  The crucial observation and one of
the main contributions of the present paper, however, is that this
results in a significant decrease in sensitivity. {\em The adaptation
law is robust with respect to perturbations in its parameters.}

\begin{figure}[h]
\hspace*{\fill}
\psfig{file=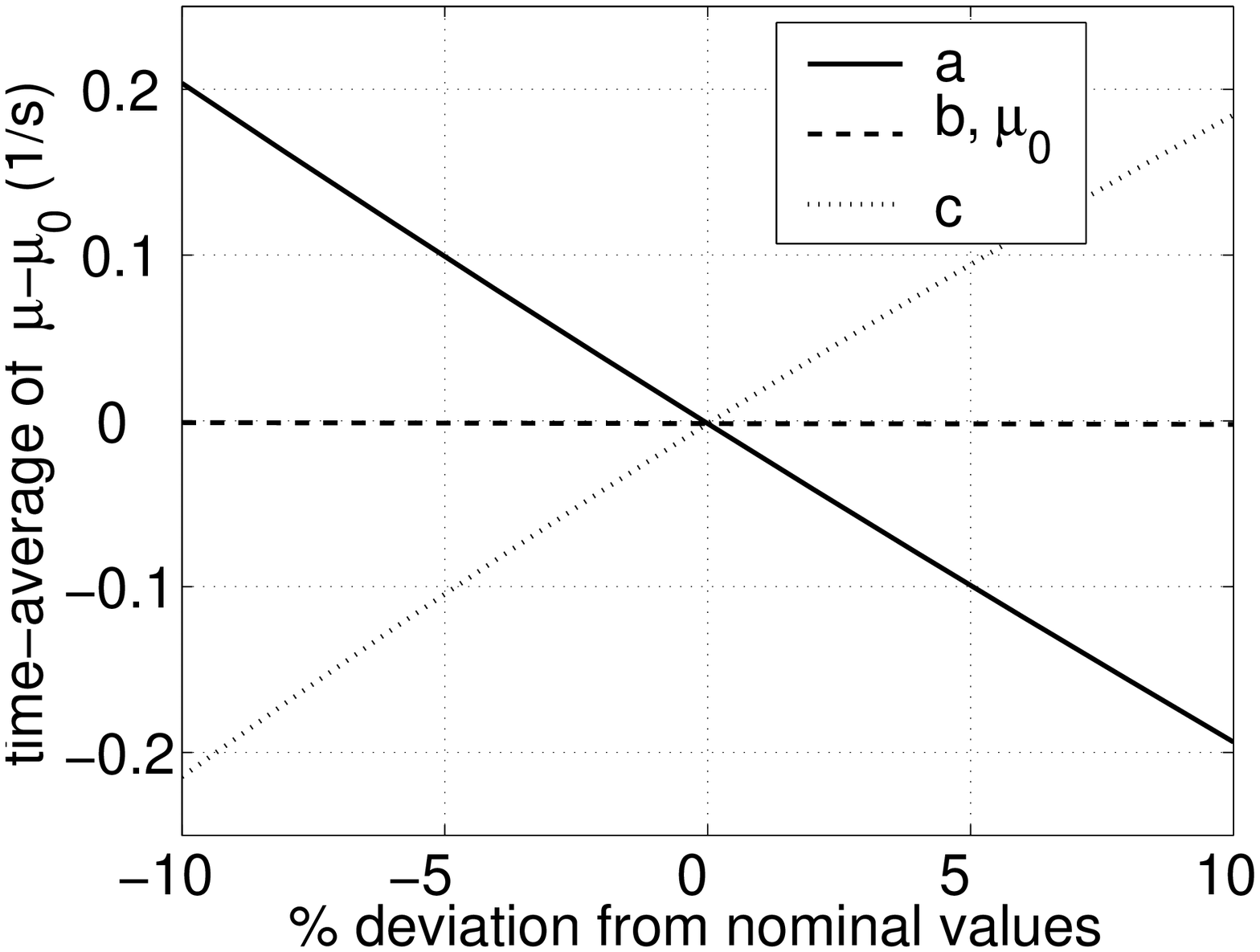,height=5cm}
\hspace*{\fill}
\psfig{file=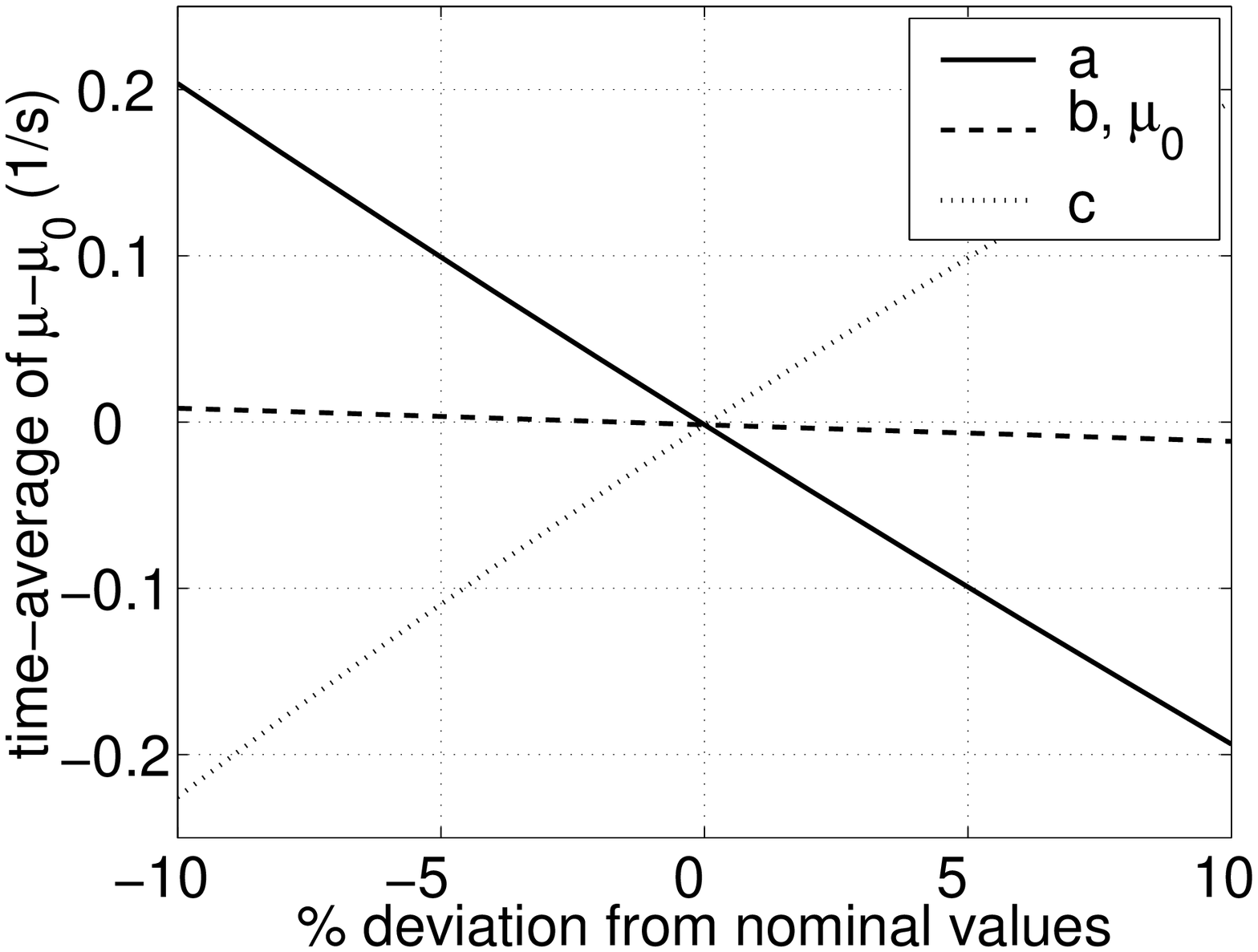,height=5cm}
\hspace*{\fill}
\caption{\label{f3} Robustness of tuning with respect to parameter
perturbations in the adaptation law.  Simulation results for
Eq.~(\ref{e:ni}) with $\dot{\mu}=0.001(-ax-b\mu+c)$ and $u(t)$ as in
Fig.~\ref{f2}.  The plots show the average value of~$\mu-\mu_0$ in
periodic regime for different values of the parameters~$a$, $b$, $c$
and~$\mu_0$.  The nominal values of the parameters satisfy the
compatibility condition~$a(20~\hertz+60~\hertz)/2+b\mu_0=c$ and are
given by $a=1~\seconds^{-1}$, $b=0.01~\seconds^{-1}$,
$c=40.1~\seconds^{-2}$, $1/\mu_0=100~\miliseconds$ (left) and
$a=1~\seconds^{-1}$, $b=0.01~\seconds^{-1}$, $c=42~\seconds^{-2}$,
$1/\mu_0=5~\miliseconds$ (right).}
\end{figure}

In order to illustrate this significant increase in robustness, let us
first summarize the well-known~\cite{SeLeReTa:neuron00} sensitivity
properties of neural integration.  Experiments suggest that the
actual time constant obtained in a tuned neural integrator circuit is
typically greater than~$10~\seconds$; that is, $|\mu-\mu_0|\leq
0.1~\seconds^{-1}$.  This requires for the fine-tuning of $\mu$ a
relative precision~$\Delta\mu/\mu$ ranging from~$1/100$ to~$1/2000$,
depending on whether the intrinsic time constant~$1/\mu_0$ equals
$100~\miliseconds$ or $5~\miliseconds$ (typical values suggested in
the literature).  The required precision for~$\mu$ should be
contrasted with the required precision for the parameters of the
adaptation law proposed in the present paper.
The simulations of
Fig.~\ref{f3} show that, in order to have $|\mu-\mu_0|\leq
0.1~\seconds^{-1}$ as observed in experiments, the parameters of the
adaptation law need to be tuned with a precision of $1/20$,
independent of the intrinsic time constant~$1/\mu_0$.
Comparing this with the originally required precision for the synaptic
feedback strength~$\mu$, we conclude that {\em the proposed adaptation
mechanism could improve the robustness of neural integration with a
factor ranging from~$5$ to~$100$}.

We have studied a simple model for neural integration and proposed a
class of feedback adaptation rules that could explain the
experimentally observed robustness of neural integration with respect
to perturbations.  The analysis tools that we have introduced extend
to the study of fine-tuning involved in other systems such as hair
cell oscillations in the ear~\cite{MoSoMu:arxiv}.
\begin{figure}[h]
\hspace*{\fill}
\psfig{file=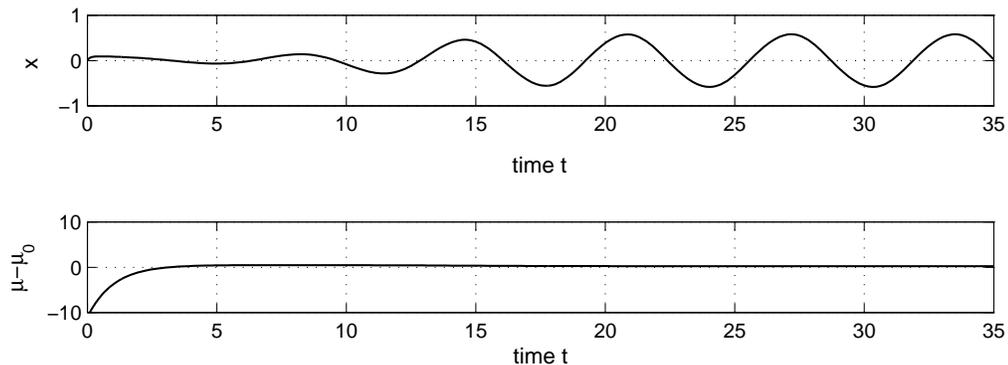,height=5cm}
\hspace*{\fill}
\caption{\label{f4} Tuning of a nonlinear oscillator.  Simulation of
equations $\ddot{x}+(\mu_0-\mu)\dot{x}+\lambda\dot{x}^3+\omega^2x=0$
and $\dot\mu=f(r)-g(\mu)$ with $\mu_0=1$, $\lambda=1$, $\omega=1$ and
$f(r)-g(\mu)=1/(1+r^2)-\mu+1/2$. The variable~$r$ is determined by
$r^2=x^2+(\dot{x}/\omega)^2$.}
\end{figure}
Consider the
nonlinear oscillator equation
$\ddot{x}+(\mu_0-\mu)\dot{x}+\lambda\dot{x}^3+\omega^2x=u(t)$, which
captures some of the essential features of hair cell
oscillations~\cite{Ju:01}.  Inspired by our previous analysis, we
consider a feedback adaptation law for the parameter~$\mu$ of the form
$\dot\mu=f(r)-g(\mu)$, with $r$ a positive variable characterizing the
magnitude of oscillations and related to~$x$ and~$\dot{x}$ via the
expression $r^2=x^2+(\dot{x}/\omega)^2$.  Fig.~\ref{f4} shows that, in
the absence of the stimulus~$u(t)$, this type of adaptation law is
indeed able to bring and keep the bifurcation parameter close to its
critical value, resulting in the spontaneous generation of
oscillations.
\section*{Acknowledgments}

Luc Moreau is a Postdoctoral Fellow of the Fund for Scientific Research -
Flanders.  Work done while
a recipient of an Honorary Fellowship of the Belgian
American Educational Foundation, while visiting the Princeton University
Mechanical and Aerospace Engineering Department.
This paper presents research results of the Belgian Programme on
Inter-University Poles of Attraction, initiated by the Belgian State,
Prime Minister's Office for Science, Technology and Culture. The
scientific responsibility rests with its authors.

Eduardo Sontag is supported in part by USAF Grant F49620-01-1-0063 and by NSF Grant
CCR-0206789.


\end{document}